\documentclass{amsart}

\usepackage{amssymb}
\usepackage{amsmath}
\newtheorem {Theorem}   {Theorem} 
\numberwithin{Theorem}{section}

\newtheorem {lemma}[Theorem]{Lemma}

\newtheorem{definition}[Theorem]{Definition}
\theoremstyle{remark}

\newtheorem{example}[Theorem]{Example}

\newtheorem{theorem}{Theorem}

\newtheorem{corollary}[theorem]{Corollary}

\begin{document}

\author{Alex Clark}
\title{Exponents and Almost Periodic Orbits}
\address{Department of Mathematics, University of North Texas, 
Denton, TX 76203-1430}
\email{alexc@unt.edu}
\date{Submitted April 1998; Revised July 1999}

\begin{abstract}We introduce the group of exponents of a map of the reals into a metric
space and give conditions under which this group embeds in the first \v
{C}ech cohomology group of the closure of the image of the map. We show that
this group generalizes the subgroup of the reals generated by the
Fourier-Bohr exponents of an almost periodic orbit and that any minimal
almost periodic flow in a complete metric space is determined up to
(topological) equivalence by this group. We also develop a way of
associating groups with any self-homeomorphism of a metric space that
generalizes the rotation number of an orientation-preserving homeomorphism
of the circle with irrational rotation number.
\vspace{.2in}
\end{abstract}

\maketitle

\section{Introduction}

In this paper we shall introduce the group of exponents of a map of $\mathbb{R}$
into a space $X$. While this group is defined for any such map, it is most
natural to consider in the case that $X$ is metric, and we will assume in
the following that all topological spaces under consideration are endowed
with a metric. We reveal the topological significance of this group by
showing that under suitable conditions it embeds (sometimes properly) in the
first \v {C}ech cohomology group of the closure of the image of the map.
This group arises naturally from a study of almost periodic orbits, and the
results of \cite{LF} together with a thorough examination of the properties of
this group allow us to give an intrinsic classification of all the minimal
almost periodic flows occurring in complete metric spaces, see Theorem \ref
{classif}. This represents progress in the program of the classification of
the limit sets of integral flows initiated by G.D. Birkhoff in 1912, see,
e.g., \cite{AM}, Chapter 6. More generally, we show that the exponent group for
the orbit of a compact minimal set determines the maximal (in a sense made
precise) semiconjugate almost periodic flow, and so this group gives an idea
of the flow shape of the minimal set.

We also calculate this exponent group for the orbits of the suspensions of
orientation--preserving homeomorphisms $S^1\rightarrow S^1$ with irrational
rotation number. This demonstrates that if we associate the exponent group
of the orbits of the suspension of any homeomorphism $X\rightarrow X$ then
we obtain groups that in some sense generalize the rotation number.

\section{The Exponent Group and its Origin}

We begin by introducing the background material on almost periodic functions
which naturally leads to the general exponent group. Throughout maps are
assumed to be continuous. By a \emph{flow} on $X$ we mean a continuous group action
of $\left( \mathbb{R},+\right) $ on $X$, and $\pi :\mathbb{R\rightarrow }\left( 
\mathbb{R}/\mathbb{Z}\right) =S^1$ denotes the standard quotient map. We shall use
the following terminology.

\begin{definition}
$S\subset \mathbb{R}$ \emph{ is } relatively dense \emph{ if there is a } $\lambda>0$ 
\emph{ such that for any } $r \in \mathbb{R}$, $\left( r,r+\lambda \right) 
\cap S\neq \emptyset .$
\end{definition}

\begin{definition} $\tau \in \mathbb{R}$ \emph{ is an } $\varepsilon$ -- almost period 
\emph{ of }  $f:\mathbb{R} \rightarrow  X$ \emph{ if }
\[
\sup_{t\in \mathbb{R}}\left\{ d\left( f\left( t\right) ,f\left( t+\tau
\right) \right)\right\} \leq \varepsilon  .
\]
\end{definition}
There is then the following definition due to H. Bohr. 
\begin{definition}
\emph{A map} $f:\mathbb{R} \rightarrow X$ \emph{ is } almost periodic  
\emph{ if for every } $\varepsilon >0 $\emph{ the set }
\[
S_{\varepsilon }=\left\{ \tau \in \mathbb{R}\mid \tau \text{\emph{ is an }}
\varepsilon -\text{\emph{almost period of }}f\right\} 
\]
\emph{is relatively dense.}
\end{definition}

	Similarly, the orbit $\phi_x$ of the flow  $\phi :\mathbb{R\times }X\rightarrow X$ 
is termed \emph{almost periodic} when it is an almost periodic function 
 $\mathbb{R} \rightarrow X .$
	
\begin{definition}
\emph{For a map }$f:\mathbb{R}\rightarrow X$\emph{, we define a sequence of
reals\ to be an }$f$--sequence \emph{if the sequence }$\{f\left( t_i\right)
\}_{i\in \mathbb{N}}$\emph{\ converges in }$X$\emph{\ and }$\{t_i\}$\emph{\ is
defined to be an }$f\left( 0\right) $--sequence\ \emph{if }$\lim_i\left\{
f\left( t_i\right) \right\} =f\left( 0\right) $\emph{.}
\end{definition}

We now quote known results translated into our terminology. One has a choice
between defining exponents $\text{mod}2\pi \,$ or $\text{mod}1$, and we
always choose the $\text{mod}1$ formulation (i.e., we use the map $\pi $ as
above instead of the map $t\mapsto \exp \left( it\right) $).

\begin{theorem}
\label{exp1} Given an almost periodic function $f:\mathbb{R}\rightarrow X$
, there exists a countable subgroup $\mathfrak{M}_f$ of $\left( 
\mathbb{R},+\right) $ such that:  
\[
\left[ \{t_i\}\subset \mathbb{R}\text{ is an }f\text{--sequence}
\right] \Leftrightarrow \left[ \{\pi \left( \lambda t_i\right) \}_i\text{
converges in } S^1 \text{ for all }\lambda \in \mathfrak{M}_f\right] ,
\]
and if $X$ is a Banach space, $\mathfrak{M}_f$ may be
taken to be the subgroup of $\mathbb{R}$ generated by the Fourier-Bohr
exponents of $f$, (see, e.g., \cite{LZ} Chapter 3, Theorems 3-4).
\end{theorem}

\begin{theorem}
\label{exp2} For two almost periodic functions $f:\mathbb{R}\rightarrow X$
 and $g:\mathbb{R}\rightarrow Y$, we have: $\left[ \mathfrak{M}_f=
\mathfrak{M}_g\right] \Leftrightarrow \left[ \{f\left( 0\right) \text{--
sequences}\}=\{g\left( 0\right) \text{-sequences}\}\right] ;$
and $\mathfrak{M}_f$ is uniquely determined by $f,$ (see, e.g.,
\cite{LZ}, Chapter 2, Section 2.3).
\end{theorem}

\begin{definition}
\label{apexp}\emph{For an almost periodic function }$f:\mathbb{R}\rightarrow X \text{, }
\mathfrak{M}_f$\emph{\ denotes the countable subgroup of }$\mathbb{R}$
\emph{\ uniquely determined by }$f\,$\emph{\ as above.}
\end{definition}

We now introduce our generalization of this group.

\begin{definition}
\label{exp}\emph{For a map }$f:\mathbb{R}\rightarrow X$\emph{, we define the }
group of exponents\emph{, denoted }$\mathcal{E}_f$\emph{, to be} 
\[
\{\alpha \in \mathbb{R\mid }\{\pi \left( \alpha t_i\right) \}\text{\emph{
converges in} }S^1\text{\emph{ for all } }f\text{-\emph{sequences} }\{t_i\}\}
\text{.}
\]
\end{definition}

\begin{lemma}
\label{subgroup}$\mathcal{E}_f$ is a subgroup of $\left( \mathbb{R},+\right) $.
\end{lemma}

\noindent \textbf{Proof:} We know $0\in \mathcal{E}_f\neq \emptyset $. Suppose that $\alpha ,\beta \in 
\mathcal{E}_f$ and that $\{t_i\}$ is an $f$--sequence. Then $\{\pi \left(
\left( \alpha -\beta \right) t_i\right) \}=\{\pi \left( \alpha t_i\right)
-\pi \left( \beta t_i\right) \}$. This sequence converges in $S^1$ since $
-\, $ is continuous on $S^1\times S^1$.\hfill$\square $

\begin{lemma}
\label{expeq}If $f\,$ is an almost periodic map into a complete metric
space, then $\mathfrak{M}_f=\mathcal{E}_f$.
\end{lemma}

\noindent \textbf{Proof:} Let $f$ be such a function. Given $\alpha \in 
\mathfrak{M}_f$, $\{\pi \left( \alpha t_i\right) \}$ converges in $S^1$ for all $
f$--sequences $\{t_i\}$, and so $\mathfrak{M}_f\subset \mathcal{E}_f$. Suppose
then that $\alpha \in \mathcal{E}_f-\mathfrak{M}_f$. Let $\mathfrak{M}_f^{\prime }$
be the subgroup of $\left( \mathbb{R},+\right) $ generated by  
$\mathfrak{M}_f\cup \{\alpha\} .$ Then
given any $f$-sequence $\left\{ t_i\right\} $, $\{\pi \left( \lambda
t_i\right) \}$ converges in $S^1$ for all $\lambda \in \mathfrak{M}_f^{\prime }$
. And if for a given sequence $\left\{ t_i\right\} $ of real numbers we have
that for each $\lambda \in \mathfrak{M}_f^{\prime }$ \ the sequence $\{\pi
\left( \lambda t_i\right) \}$ converges in $S^1$, then for each $\lambda \in 
\mathfrak{M}_f\,$\ the sequence $\{\pi \left( \lambda t_i\right) \}$ converges
in $S^1$ since $\mathfrak{M}_f\subset \mathfrak{M}_f^{\prime }$, and so $\left\{
t_i\right\} $ is an $f-$sequence. Therefore $\mathfrak{M}_f^{\prime }$ meets the
conditions that uniquely determine $\mathfrak{M}_f$,\textbf{\ }and so \textbf{\ 
}$\mathfrak{M}_f^{\prime }=\mathfrak{M}_f$, contradicting the choice of $\alpha $.
Thus, no such $\alpha $ can exist and $\mathcal{E}_f=\mathfrak{M}_f$.\hfill$\square 
$

Now we shall explore some properties of this more general exponent group. Here $\left[
M;S^1\right] $ denotes the group of homotopy classes of maps $M \rightarrow S^1$, where
for convenience we use the group structure induced by point-wise addition of maps.

\begin{definition}
\emph{We define a map }$f:\mathbb{R}\rightarrow X$\emph{\ to be }stable\emph{\
if there is an unbounded }$f$\emph{-sequence }$\{t_i\}$\emph{.}
\end{definition}

\begin{theorem}
\label{embedexp}If $f:\mathbb{R}\rightarrow X$ is a stable map and $M=\overline{
f\left( \mathbb{R}\right) }$, the map $\iota :\mathcal{E}_f\rightarrow \left[
M;S^1\right] $ given by 
\[
\alpha \mapsto \left[ f_\alpha \right] ,\text{ where }f_\alpha :M\rightarrow
S^1;\;f_\alpha \left( \{\lim f\left( t_i\right) \}\right) =\lim \{\pi \left(
\alpha t_i\right) \}\text{ }
\]
\hspace{.58in}and $\left[ f_\alpha \right] $ denotes the homotopy class of $
f_\alpha $, \
\newline is an embedding. 
\end{theorem}

\noindent \textbf{Proof:}

\begin{enumerate}
\item  $\iota $ is well-defined:
\end{enumerate}

Let $\alpha \in \mathcal{E}_f$ be given. If $\lim \{f\left( t_i\right)
\}=\lim \{f\left( t_i^{\prime }\right) \}=x\in M$ are two representations of
a point of $M$, then $\lim \{\pi \left( \alpha t_i\right) \}=\xi $ and $\lim
\{\pi \left( \alpha t_i^{\prime }\right) \}=\xi ^{\prime }$ both exist by
the definition of $\mathcal{E}_f$. Then with $s_{2i}=t_i$ and $
s_{2i-1}=t_i^{\prime }$, $\lim \{f\left( s_i\right) \}=x\Rightarrow \lim
\{\pi \left( \alpha s_i\right) \}$ exists, which is only possible if $\xi
=\xi ^{\prime }$. Thus $f_\alpha $ is a well-defined function $M\rightarrow
S^1$. To see that $f_\alpha $ is continuous, consider a convergent sequence $
\{\lim\limits_i\{f\left( t_i^j\right) \}\}_j=\{x^j\}_j\rightarrow x$. Then 
for each $j\in \mathbb{N}$ we may choose $i_j$ so that: 
\[
d\left( f\left( t_{i_j}^j\right) ,x^j\right) <\frac 1j\text{ and }
d_{S^1}\left( f_\alpha \left( x^j\right) ,\pi \left( \alpha t_{i_j}^j\right)
\right) <\frac 1j\text{.} 
\]
Then $\lim\limits_j\{f\left( t_{i_j}^j\right) \}_j=x$ and 
\[
f_\alpha \left( x\right) \stackrel{\text{def}}{=}\lim\limits_j\{\pi \left( \alpha
t_{i_j}^j\right) \}=\lim\limits_j\{f_\alpha \left( x^j\right) \}\text{, } 
\]
demonstrating that $f_\alpha $ is continuous.

\begin{enumerate}
\item[2.]  $\iota $ is a homomorphism:
\end{enumerate}

Let $\alpha ,\beta \in \mathcal{E}_f$. Then
for $x=\lim \{f\left( t_i\right) \}\in M$ 
\[
f_{\alpha +\beta }(x)=\lim \{\pi \left( \left( \alpha +\beta \right)
t_i\right) \}=\lim \{\pi \left( \alpha t_i\right) \}+\lim \{\pi \left( \beta
t_i\right) \}=f_\alpha (x)+f_\beta (x)\text{.} 
\]

\begin{enumerate}
\item[3.]  $\ker \iota =\{0\}$:
\end{enumerate}

Suppose $0\neq \alpha \in \mathcal{E}_f$ and $[f_\alpha ]=\left[ \text{constant map}\right] .$
Since $\pi :\mathbb{R\rightarrow }S^1$ is a fibration
and since we can lift the constant map, we can lift $f_\alpha $ with a map $
g:(M,f(0))\rightarrow (\mathbb{R},0)$ making the following diagram commute:

\begin{equation*}
\begin{array}{ccc}
&  & \mathbb{R} \\ 
& \overset{g}{\nearrow } & \downarrow \pi  \\ 
M & \underset{f_{\alpha }}{\longrightarrow } & S^{1}
\end{array}
.
\end{equation*}

This leads to the following commutative diagram: 
\[
\begin{array}{ccccccc}
\mathbb{R} & \stackrel{\frac 1\alpha }{\rightarrow } & \mathbb{R} & \stackrel{f}{
\rightarrow } & M & \stackrel{g}{\rightarrow } & \mathbb{R} \\ 
\downarrow ^\pi &  &  &  &  & 
\begin{array}{cc}
\searrow^{f_\alpha}
\end{array}
& \downarrow ^\pi \\ 
S^1 &  &  & \stackrel{id}{\longrightarrow } &  &  & S^1
\end{array}
, 
\]
where $gf\frac 1\alpha :(\mathbb{R},0)\mathbb{\rightarrow }(\mathbb{R},0);$ $t\mapsto
g\circ f\left( \frac t\alpha \right) $ . From this it follows that $gf\frac
1\alpha $ is a lift of $id_{S^1}$. Such a lift $(\mathbb{R},0)\mathbb{\rightarrow }
(\mathbb{R},0)$ is uniquely determined \cite{S}, 2.2 Lemma 4 (specifying that $
0\mapsto 0$ makes the lift unique), and $id_{\mathbb{R}}$ also provides such a
lift. Thus $gf\frac 1\alpha =id_{\mathbb{R}}$. Composing both sides of this
equality on the right with multiplication by $\alpha $, we obtain: $gf=$
multiplication by $\alpha $. Since $f$ is stable, there is an unbounded
sequence of numbers $\{t_i\}$ such that $\lim \{f\left( t_i\right) \}=x$.
And so 
\[
g\left( x\right) =g\left( \lim \{f\left( t_i\right) \}\right) =\lim \{g\circ
f(t_i)\}=\lim \{\alpha t_i\}\emph{,} 
\]
which is not well-defined since $\{t_i\}$ and hence $\{\alpha t_i\}$ is an
unbounded sequence that does not converge. This contradicts the continuity
of $g$.\hfill$\square $

Notice that when $M$ is compact $f$ is automatically stable and $\mathcal{E}
_f$ is countable since $\left[ M;S^1\right] $ is then countable. While we
shall be principally interested in maps which are orbits of a flow, we give
a more general example first.

\begin{example}
For $t\in \mathbb{R}$, $\left\lfloor t\right\rfloor \stackrel{\text{def}
}{=}$(the greatest integer less than $t$) and $\left\rfloor
t\right\lfloor \stackrel{\text{def}}{=}t-\left\lfloor t\right\rfloor .$ Then
we define the map $g:\mathbb{R\rightarrow R}^2;$

\begin{equation*}
t\mapsto \left\{ 
\begin{array}{ll}
\left( \rfloor t\lfloor \cdot \frac{1}{2^{\lfloor t\rfloor +1}}+(1-\rfloor
t\lfloor )\cdot \frac{1}{2^{\lfloor t\rfloor }},\rfloor t\lfloor \right)  & 
\text{if }\lfloor t\rfloor \text{ is odd} \\ 
\left( \rfloor t\lfloor \cdot \frac{1}{2^{\lfloor t\rfloor +1}}+(1-\rfloor
t\lfloor )\cdot \frac{1}{2^{\lfloor t \rfloor }},1-\rfloor t\lfloor \right) 
& \text{if }\lfloor t\rfloor \text{ is even}
\end{array}
\right. 
\end{equation*}

\noindent and $X\stackrel{\text{def}}{=}g\left( \mathbb{R}\right) \cup \{\left( 0,\frac
12\right) \}$, which is homeomorphic to 
\[
\left\{ \left( x,\sin \left(
1/x\right) \right) \mid x\in (0,\infty )\right\} \cup \left\{ \left(
0,0\right) \right\} \subset \mathbb{R}^2.
\] 

\noindent Then $f:\mathbb{R\rightarrow }
X;t\mapsto g\left( t\right) $ is stable since $\{f\left( n+\frac
12\right) \}_{n\in \mathbb{N}}\rightarrow \left( 0,\frac 12\right) $ . And
if $\{f\left( t_i\right) \}\rightarrow x\in X$, we have
two cases:
\end{example}

\begin{enumerate}
\item  $x=f\left( t\right) $ for some $t\in \mathbb{R}$.
\end{enumerate}

In this case we have the neighborhood $U\stackrel{\text{def}}{=}\{\left( w,y\right)
\in $ $X\mid \frac 1{2^{^{\left\lfloor t\right\rfloor +1}}}\leq w\leq \frac
1{2^{\left\lfloor t\right\rfloor -1}}\}$ of $f\left( t\right) $ and $f$ maps 
$\left[ \left\lfloor t\right\rfloor -1,\left\lfloor t\right\rfloor +1\right] 
$ homeomorphically onto $U$. Thus, we must have $\left\{ t_i\right\}
\rightarrow t$ and so $\left\{ \pi \left( 1\cdot t_i\right) \right\}
\rightarrow \left\{ \pi \left( 1\cdot t\right) \right\} $ in $S^1$.

\begin{enumerate}
\item[2.]  $x=\left( 0,\frac 12\right) $.
\end{enumerate}

In this case, with $k_i\stackrel{\text{def}}{=}\left\lfloor t_i\right\rfloor $ and $
\delta _i$ $\stackrel{\text{def}}{=}\left\rfloor t_i\right\lfloor $ , we have that $
\left\{ \delta _i\right\} \rightarrow \frac 12$, and so $\left\{ \pi \left(
1\cdot t_i\right) \right\} =\left\{ \pi \left( \delta _i\right) \right\}
\rightarrow \pi \left( \frac 12 \right)$ in $S^1$.

Thus, we must have $1$ and hence $\mathbb{Z}\subset \mathcal{E}_f$. Let $\alpha
=\frac pq\in \mathbb{Q}-\mathbb{Z}$. We then construct the $f$-sequence $\{t_i\}$,

\begin{equation*}
t_{i}\stackrel{\text{def}}{=}\left\{ 
\begin{array}{ll}
qi+\frac{1}{2} & \text{if }i\text{ is odd} \\ 
qi+\frac{3}{2} & \text{if }i\text{ is even}
\end{array}
.\right. 
\end{equation*}
 
Then
\begin{equation*}
\pi \left( \alpha \cdot t_{i}\right) =\left\{ 
\begin{array}{ll}
\pi \left( \frac{p}{2q}\right)  & \text{if }i\text{ is odd} \\ 
\pi \left( \frac{p}{2q}+\frac{p}{q}\right)  & \text{if }i\text{ is even}
\end{array}
\right.
\end{equation*}

\noindent and $\pi \left( \frac p{2q}\right) \neq \pi \left( \frac p{2q}+\frac
pq\right) $ in $S^1$ since $\frac pq\in \mathbb{Q}-\mathbb{Z}$, and so $\left\{ \pi
\left( \alpha \cdot t_i\right) \right\} $ does not converge in $S^1\,$ and $
\alpha \notin \mathcal{E}_f$. Suppose then that $\alpha \in \mathbb{R}-\mathbb{Q}$
. By Kronecker's Theorem (see, e.g., \cite{LZ},
Chapter 3.1) there are then two sequences of integers $\{k_i\}$
and $\left\{ \ell _i\right\} $ satisfying

\begin{enumerate}
\item  $\{\pi \left( \alpha \cdot k_i\right) \}\rightarrow \pi \left(0\right)$ 
in $S^1\,$ and 
$\{\pi \left( \alpha \cdot \ell _i\right) \}\rightarrow \pi \left(\frac 13\right)$ in $S^1$

\item  $k_i,\ell _i>i$ for all $i\in \mathbb{N}$.
\end{enumerate}

Then with
\begin{equation*}
t_{i}\overset{def}{=}\left\{ 
\begin{array}{ll}
k_{i}+\frac{1}{2} & \text{if }i\text{ is odd} \\ 
\ell _{i}+\frac{1}{2} & \text{if }i\text{ is even}
\end{array}
\right. ,
\end{equation*}
we have that $ \left\{ f\left( t_i\right) \right\} \rightarrow \left( 0,\frac 12\right) $
while $\left\{ \pi \left( \alpha \cdot t_i\right) \right\} $ does not
converge in $S^1$. Thus, $\mathbb{Z}=\mathcal{E}_f$ and $\left[ X;S^1\right] \,$
contains a copy of $\mathbb{Z}$.

It is perhaps of interest to calculate $\breve{H}^1\left( X\right) \cong
\left[ X;S^1\right] $: we can find a cofinal sequence of nerves consisting
of a line attached to the common point of a bouquet of countably infinitely many
circles in the weak topology, with each circle collapsing to the line in all
refinements after a certain point in the sequence. And so $\breve{H}^1\left(
X\right) $ is the direct limit of a sequence of groups isomorphic with $
Hom\left( \oplus _{i=1}^\infty \mathbb{Z},\mathbb{Z}\right) \cong
\prod_{i=1}^\infty \mathbb{Z}$ with the bonding maps sending a finite number of 
$\mathbb{Z}$ factors to $0$ in the group farther along in the sequence, and so $
\breve{H}^1\left( X\right) $ is isomorphic to the group of sequences of
integers with sequences identified which are eventually the same.

\subsection{$\mathfrak{M}_f$ Determines the Equivalence Class of $\phi $}

In the sequel we shall assume all spaces are complete in their endowed
metric unless otherwise stated. The following basic result does not seem to
be proven anywhere in the literature, but it seems to be implicitly assumed in 
\cite{LZ}. Cartwright \cite{C} does observe that,
``A flow is uniquely determined by its coefficients and its exponents.''
Notice that the ``coefficients'' need not be brought into the picture when
determining the equivalence class of the corresponding almost periodic flow.
We shall use the following terminology and notation.

\begin{definition}
\label{semicon}\emph{The flow }$\psi $\emph{\ on }$Y$\emph{\ is }
semiconjugate\emph{\ to the flow }$\phi $\emph{\ on }$X$\emph{\ if there
exists a surjective map }$h:X\rightarrow Y\,$\emph{\ such that the following
diagram commutes} 
\[
\begin{array}{ccc}
\mathbb{R}\times X & \stackrel{\phi }{\longrightarrow } & X \\ 
^{\left( \alpha ,h\right) }\downarrow  &  & \downarrow ^h \\ 
\mathbb{R}\times Y & \stackrel{\psi }{\longrightarrow } & Y
\end{array}
\emph{,}
\]
\emph{where }$\alpha $\emph{\ is multiplication by }$a\in \mathbb{R-}\left\{
0\right\} $\emph{, and we write }$\alpha \times h:\phi \stackrel{sc}{\succeq 
}\psi $\emph{. And if }$h$\emph{\ is a homeomorphism and if }$a>0$\emph{, we
write }$\alpha \times h:\phi \stackrel{equiv}{\approx }\psi $\emph{\ (or }$
\phi \stackrel{equiv}{\approx }\psi $\emph{) and we say that }$\phi $\emph{\
and }$\psi $\emph{\ are }equivalent\emph{.}
\end{definition}

\begin{theorem}
\label{class}If $f=\phi _x:\mathbb{R}\rightarrow X$ is an almost periodic orbit
of a flow $\phi :\mathbb{R\times }X\rightarrow X$ on $(X,d)$ and $g=\psi _y:
\mathbb{R}\rightarrow Y$ is an almost periodic orbit of a flow $\psi :\mathbb{
R\times }Y\rightarrow Y$ on $(Y,d^{\prime })$ and $\mathfrak{M}_f$ $=$ $\mathfrak{M}
_g$, then $id_{\mathbb{R}}\times h:\phi \mid _{\mathbb{R\times }\overline{f(\mathbb{R)
}}}\stackrel{equiv}{\approx }$ $\psi \mid _{\mathbb{R\times }\overline{g(\mathbb{R)
}}}$, where $h$ is a topological isomorphism of $\overline{f(\mathbb{R)}}$ and $
\overline{g(\mathbb{R)}}$ with the group structures inherited from the orbits $f
$ and $g$.
\end{theorem}

\noindent \textbf{Proof}: The group operations on $($ $\overline{f(\mathbb{R)}}
,+$ $)$ and $($ $\overline{g(\mathbb{R)}},+)$ discovered by Nemytski \cite{NS}, V, Thm
8.16, are given by $\lim \{f(t_i)\}+\lim \{f(t_i^{\prime })\}=\lim
\{f(t_i+t_i^{\prime })\}$ and $\lim \{g(t_i)\}+\lim \{g(t_i^{\prime
})\}=\lim \{g(t_i+t_i^{\prime })\}$ respectively . The continuity of $\phi $
then yields that $\phi \left( t,\lim \{f(t_i)\}\right) =\phi \left( \lim
\{\left( t,f\left( t_i\right) \right) \}\right) =\lim \{\phi \left(
t,f\left( t_i\right) \right) \}=\lim \{f(t+t_i)\}$. Define $h:$ $\overline{f(
\mathbb{R)}}\rightarrow $ $\overline{g(\mathbb{R)}}$ by $\lim \{f(t_i)\}\mapsto
\lim \{g(t_i)\}$. Since $\mathfrak{M}_f$ = $\mathfrak{M}_g$, $\left[ \lim \{f(t_i)\}
\text{ exists}\Leftrightarrow \lim \{g(t_i)\}\text{ exists}\right] $.
Therefore, $h$ will be well-defined if we can show that when $a=\lim
\{f(t_i)\}$ and $a=\lim \{f(t_i^{\prime })\},$ then $\lim \{g(t_i)\}=\lim
\{g(t_i^{\prime })\}$. If $a=\lim \{f(t_i)\}$ and $a=\lim \{f(t_i^{\prime
})\},$ $f(0)=a-a=\lim \{f(t_i-t_i^{\prime })\}$, which implies that $
\{t_i-t_i^{\prime }\}$ is an $f\left( 0\right) $-sequence and hence a $
g\left( 0\right) $-sequence, see Theorem \ref{exp2}. From this it follows that 
\[
\lim \{g(t_i)\}-\lim \{g(t_i^{\prime })\}=\lim \{g(t_i-t_i^{\prime })\}=g(0)
\]
and hence that $\lim \{g(t_i)\}=\lim \{g(t_i^{\prime })\}$. Thus, $h$ is
well-defined and a similar argument shows that $k:\overline{g(\mathbb{R)}}
\rightarrow \overline{f(\mathbb{R)}};$\ $\lim \{g(t_i)\}\mapsto \lim \{f(t_i)\}$
is well-defined, and $h\circ k=id_{\overline{g(\mathbb{R)}}}$ and $k\circ h=id_{
\overline{f(\mathbb{R)}}}$, demonstrating that $h$ is a bijection.

If 
\begin{equation*}
a=\lim \left\{ a_{j}\right\} =\lim_{j}\left\{ \lim_{i}\left\{ f\left(
t_{i}^{j}\right) \right\} \right\} 
\end{equation*}
in $\overline{f(\mathbb{R)}}$, then 
\begin{equation*}
a=\lim_{j}\left\{ f\left( t_{i_{j}}^{j}\right) \right\} 
\end{equation*}
for $i_j$ satisfying $d\left( f\left(
t_n^j\right) ,a_j\right) <\frac 1j$ and $d^{\prime }\left( g\left(
t_n^j\right) ,h\left( a_j\right) \right) <\frac 1j$ for all $n\geq $ $i_j$.
Such $i_j$ exist since

\begin{equation*}
a_{j}=\lim_{i}\left\{ f\left( t_{i}^{j}\right) \right\}
\text{ and }h\left(
a_{j}\right) =\lim_{i}\left\{ g\left( t_{i}^{j}\right) \right\} .
\end{equation*}
Then 
\begin{equation*}
h\left( a\right) =\lim_{j}\left\{ g\left( t_{i_{j}}^{j}\right) \right\}
=\lim_{j}\left\{ \lim_{i}\left\{ g\left( t_{i}^{j}\right) \right\} \right\}
=\lim_{j}\left\{ h\left( a_{j}\right) \right\} ,
\end{equation*} and so $h$ is continuous and therefore a homeomorphism. And since 
\begin{eqnarray*}
h\left( \lim \{f(t_i)\}+\lim \{f(t_i^{\prime })\}\right) &=&h\left( \lim
\{f(t_i+t_i^{\prime })\}\right) =\lim \{g(t_i+t_i^{\prime })\}= \\
\lim \{g(t_i)\}+\lim \{g(t_i^{\prime })\} &=&h\left( \lim \{f(t_i)\}\right)
+h\left( \lim \{f(t_i^{\prime })\}\right) ,
\end{eqnarray*}
we see that $h$ is a homomorphism and thus a topological isomorphism. Also: 
\begin{eqnarray*}
h\left( \phi \left( t,\lim \{f(t_i)\}\right) \right) &=&h\left( \lim
\{f\left( t+t_i\right) \}\right) =\lim \{g\left( t+t_i\right) \}= \\
\psi \left( t,\lim \{g(t_i)\}\right) &=&\psi \left( t,h\left( \lim
\{f(t_i)\}\right) \right) ,
\end{eqnarray*}
from which it follows that $id_{\mathbb{R}}\times h$ provides the desired
equivalence.\hfill$\square $

Notice that the examples $X$ = one trajectory of an irrational flow on the
torus and $Y$ = two trajectories of the same irrational flow on the torus
show that the theorem does not hold as stated without the completeness
requirement.

\section{Properties of $\mathcal{E}_f$ in Flows}

\begin{theorem}
\label{minexp}If $f\,$ and $g$ are the orbits of $x\,$ and $y\,$
respectively of the flow $\phi $ on $X$ and if $y\in \overline{f\left( \mathbb{R
}\right) }$, then $\mathcal{E}_f\subset \mathcal{E}_g$. If both $x$ and $y$
are contained in a minimal set $M$, then $\mathcal{E}_f=\mathcal{E}_g$.
\end{theorem}

\noindent \textbf{Proof}: Let $\left\{ f\left( t_i\right) \right\}
\rightarrow y$ and suppose $\left\{ s_i\right\} \,$ is a $g$-sequence with $
\left\{ g\left( s_i\right) \right\} \rightarrow \xi $. Then 
\begin{eqnarray*}
\xi &=&\lim_i\left\{ \phi \left( s_i,y\right) \right\} =\lim_i\left\{ \phi
\left( s_i,\lim_j\left\{ \phi \left( t_j,x\right) \right\} \right) \right\}
=\lim_i\left\{ \phi \left( s_i,\phi \left( t_{j_i},x\right) \right) \right\}
\\
&=&\lim_i\left\{ \phi \left( s_i+t_{j_i},x\right) \right\} =\lim_i\left\{
f\left( s_i+t_{j_i}\right) \right\}
\end{eqnarray*}
for an appropriately chosen subsequence $\left\{ t_{j_i}\right\} $ of $
\left\{ t_j\right\} $. Then we have that both $\left\{ t_{j_i}\right\} $ and 
$\left\{ s_i+t_{j_i}\right\} $ are $f$-sequences. Let $\lambda \in \mathcal{E
}_f$. Then we have that both $\left\{ \pi \left( \lambda t_{j_i}\right)
\right\} _i$ and $\left\{ \pi \left( \lambda \left( s_i+t_{j_i}\right)
\right) \right\} _i=\left\{ \pi \left( \lambda s_i\right) +\pi \left(
\lambda t_{j_i}\right) \right\} _i$ converge. But then the sequence $\left\{
\pi \left( \lambda s_i\right) \right\} _i=\left\{ \pi \left( \lambda
s_i\right) +\pi \left( \lambda t_{j_i}\right) -\pi \left( \lambda
t_{j_i}\right) \right\} _i$ converges since $-$ is continuous on $S^1\times
S^1$. Thus, $\mathcal{E}_f\subset \mathcal{E}_g$. If both $x$ and $y$ are
contained in a minimal set $M$, then $x\in \overline{g\left( \mathbb{R}\right) }
=M$ and $y\in \overline{f\left( \mathbb{R}\right) }=M$ and so we have both $
\mathcal{E}_g\subset \mathcal{E}_f$ and $\mathcal{E}_f\subset \mathcal{E}_g$.
\hfill $\square$

We now provide an example to show that the containment can be proper.

\begin{example}
\label{spiral}Let $X=S^1\times \mathbb{R}$ be parameterized by $
\left( \theta ,r\right) $ and for non--zero real numbers $\alpha $ and $
\beta $ let $\phi \left( \alpha ,\beta \right) $ be the flow on
$X\,$ generated by the vector field 
\[
\frac{d\theta }{dt}=\alpha r+\beta \left( 1-r\right) \,\text{ and }
\frac{dr}{dt}=r\left( 1-r\right) \text{.}
\]
\end{example}

Then the orbit $g\,$ of the point $\left( \pi \left( 0\right) ,1\right) $ is
given by $t\mapsto \left( \pi \left( \alpha t\right) ,1\right) $ and the
orbit $g^{\prime }$ of the point $\left( \pi \left( 0\right) ,0\right) $ is
given by $t\mapsto \left( \pi \left( \beta t\right) ,0\right) $. And so $
\mathcal{E}_g=\left\langle \alpha \right\rangle _{\mathbb{Z}}$ and $\mathcal{E}
_{g^{\prime }}=\left\langle \beta \right\rangle _{\mathbb{Z}}$. The orbit $f$
of the point $\left( \pi \left( 0\right) ,\frac 12\right) $ is given by 
\[
t\mapsto \left( \pi \left( \alpha \ln \left( e^t+1\right) -\beta \ln \left(
e^{-t}+1\right) +\left( \beta -\alpha \right) \ln 2\right) ,\frac{e^t}{e^t+1}
\right) . 
\]
And so $\left\{ \left( \pi \left( 0\right) ,1\right) ,\left( \pi \left(
0\right) ,0\right) \right\} \subset \overline{f\left( \mathbb{R}\right) }$. If
we choose $\alpha $ and $\beta $ to be rationally independent then we have by
the above $\mathcal{E}_f\subset \mathcal{E}_g\cap \mathcal{E}_{g^{\prime
}}=\left\{ 0\right\} $, and so $\mathcal{E}_f=\left\{ 0\right\} $ and we have
proper containment in this case. We mention in passing that we can define
exponent groups $\mathcal{E}_f^{+}$ for semi-orbits $f:[0,\infty
)\rightarrow X$ in almost exactly the same way as we have defined exponents
for orbits, and when we do so, for the semi-orbit of the $f\,$ in this
example we obtain: $\mathcal{E}_f^{+}=\left\langle \alpha \right\rangle _{
\mathbb{Z}}$, the subgroup of $\left( \mathbb{R},+ \right)$ generated by $\alpha$ .

\begin{lemma}
\label{dominate}If $\alpha \times h:\phi \stackrel{sc}{\succeq }\psi $, then
for any $x\in X$ we have: 
\[
\mathcal{E}_{\phi _x}\supset \left\{ a\lambda \mid \lambda \in \mathcal{E}
_{\psi _{h\left( x\right) }}\right\} \text{.}
\]
And so $\mathcal{E}_{\phi _x}\supset \mathcal{E}_{\psi _{h\left( x\right) }}$
if $a=1$.
\end{lemma}

\noindent \textbf{Proof}: By hypothesis we have $h\circ \phi _x=\psi
_{h\left( x\right) }\circ \alpha $. Suppose $\left\{ t_i\right\} $ is a $
\phi _x$-sequence. Then $\lim\limits_i\left\{ \psi _{h\left( x\right)
}\left( at_i\right) \right\} =\lim\limits_i\left\{ h\circ \phi _x\left(
t_i\right) \right\} $ and so $\left\{ at_i\right\} $ is a $\psi _{h\left(
x\right) }$-sequence. Then for any $\lambda \in \mathcal{E}_{\psi _{h\left(
x\right) }}$ we have that $\left\{ \pi \left( a\lambda t_i\right) \right\}
_i $ converges.\hfill$\square $

\begin{corollary}
\label{equivexp}If $\left( \alpha \times h\right) :\phi \stackrel{equiv}{
\approx }\psi $, then for any $x\in X$ we have $\mathcal{E}_{\phi
_x}=\left\{ a\lambda \mid \lambda \in \mathcal{E}_{\psi _{h\left( x\right)
}}\right\} $. And so $\mathcal{E}_{\phi _x}=\mathcal{E}_{\psi _{h\left(
x\right) }}$ if $a=1$.\hfill$\square $
\end{corollary}

\section{Constructing a Flow Dual to a Given Exponent Group}

For the sake of uniformity, we shall make use of the terminology and
notation found in \cite{F}. Here $\mathbf{T}
^\kappa $ denotes the $\kappa$--fold product of $S^1$ with points $ \langle x_1,x_2,
...\rangle$, and $\Phi _{\overline{M}}^{\mathbf{\omega }}$ is the linear flow on 
$\sum_{\overline{M}}$
\[
\left( t, x \right) \mapsto \pi_{\overline{M}} \left(\omega t\right) + x,
\]
as defined in \cite{LF}.

\begin{definition}
\label{Bseq}\emph{Given a countable subgroup of the reals }$H$\emph{\ }$
=\{h_1=0,h_2,...\}$\emph{\ with a maximal independent set }$
B=\{b_i\}_{i=1}^\kappa $\emph{, we define the }$B$-sequence\emph{\ of }$H$
\emph{\ to be the direct sequence }$\{H^i,\beta _i^j\}$\emph{\ with } 
\[
H^1=\langle B\rangle _{\mathbb{Z}},H^2=\langle B\cup \{h_2\}\rangle _{\mathbb{Z}
},...,H^n=\langle H^{n-1}\cup \{h_n\}\rangle _{\mathbb{Z}},...
\]
\emph{and with }$\beta _i^j:H^i\hookrightarrow H^j$\emph{\ inclusion
(\thinspace here }$\langle S\rangle _{\mathbb{Z}}$\emph{\ denotes the subgroup
of }$H$\emph{\ generated by the set }$S$\emph{\thinspace ).}
\end{definition}

By standard results, $H$ is isomorphic to the direct limit of any $B$
-sequence of $H$. Also, $\kappa =r_0\left( H\right) $
(the torsion-free rank of $H$), and so  the cardinality of $B$ (namely, $\kappa$)
is uniquely determined by $H$ \cite{F}, III, 16.3.

\begin{lemma}
If $\{H^i,\beta _i^j\}$ is the $B$-sequence of the countable subgroup of the
reals $H=\{h_1=0,h_2,...\}$, and if $B=\{b_i\}_{i=1}^\kappa $ has
cardinality $\kappa $, then for each $i\in \mathbb{N}$ the group $H^i$ is
isomorphic to $\bigoplus_{j=1}^\kappa \langle b_j^i\rangle _{\mathbb{Z}}\cong
\bigoplus_{j=1}^\kappa \mathbb{Z}$ for some $b_j^i\in H^i$ and each bonding map 
$\beta _i^{i+1}$ can be represented by: (1) an $n\times n\,$ invertible
matrix $\left( M_i\right) ^T$ with integer entries if $\kappa =n<\infty $ or
(2) a map $\left( M_i\right) ^T\times id$, where $M_i$ is an invertible $
n_i\times n_i\,$ matrix with integer entries if $\kappa =\infty $.
\end{lemma}

\noindent \textbf{Proof:} We shall only treat the case $\kappa =\infty $,
the finite case being handled similarly. Proceeding by induction, as $
\{b_j^{i-1}\}_{j=1}^\infty \,$ is a maximal independent set of $H$, we have
that $h_i=r_1b_1^{i-1}+\cdots +r_kb_k^{i-1}$ for some $k\in \mathbb{N}$ and
rationals $r_1,...,r_k$. With $M^{\prime }\stackrel{\text{def}}{=}
\{b_j^{i-1}\}_{j=1}^k$ and $A\stackrel{\text{def}}{=}\langle M^{\prime }\cup
\{h_i\}\rangle _{\mathbb{Z}}$, we have that $M^{\prime }$ is an independent
subset of $A$ and any element of $\langle h_i\rangle _{\mathbb{Z}}$ depends on $
M^{\prime }$. Thus, $M^{\prime }\,$ is a maximal independent subset of $A$.
The finitely generated group $A$ is the direct sum of a finite number $
\left( k^{\prime }\right) $ of cyclic groups, and the invariance of $
r_0\left( A\right) =\left| M^{\prime }\right| $ implies that $k^{\prime
}=\left| M^{\prime }\right| =$ $k$:\ $A\cong \bigoplus_{j=1}^k\langle
b_j^i\rangle _{\mathbb{Z}}\cong \bigoplus_{j=1}^k\mathbb{Z}$, for some $
b_1^i,...,b_k^i$ contained in $A$. Then for $j>k\,$, let $b_j^i\stackrel{\text{def}
}{=}b_j^{i-1}$. We need to show that $\{b_j^i\}_{j=1}^\infty $ is an
independent subset of $H^i$. If we have a relation 
\[
\left( \ast \right) \text{ }n_1b_1^i+\cdots +n_kb_k^i+\cdots +n_\ell b_\ell
^i=0, 
\]
where $n_s\in \mathbb{Z}$, then for each $j\leq k$ we have $
b_j^i=r_1^jb_1^{i-1}+\cdots +r_k^jb_k^{i-1}$ for rational numbers $r_s^j$
since $h_i$ may be expressed this way and $b_j^i\in \langle M^{\prime }\cup
\{h_i\}\rangle _{\mathbb{Z}}$. Our original relation then leads to a relation
of the form $\rho _1b_1^{i-1}+\cdots +\rho
_kb_k^{i-1}+n_{k+1}b_{k+1}^{i-1}+\cdots +n_\ell b_\ell ^{i-1}=0$, where each 
$\rho _s$ is rational. Multiplying through by $\lambda \stackrel{\text{def}}{=}(
\text{lcm}\,$ of the denominators of the $\rho _s)$, we obtain: 
\[
\nu _1b_1^{i-1}+\cdots +\nu _kb_k^{i-1}+\lambda n_{k+1}b_{k+1}^{i-1}+\cdots
+\lambda n_\ell b_\ell ^{i-1}=0, 
\]
where $\nu _s\in \mathbb{Z}$. The independence of $\{b_j^{i-1}\}_{j=1}^\kappa $
then yields that 
\[
\lambda n_{k+1}=\cdots =\lambda n_\ell =0\Rightarrow n_{k+1}=\cdots =n_\ell
=0. 
\]
Our original relation $\left( *\right) $ is then of the form $
n_1b_1^i+\cdots +n_kb_k^i=0$, and since $\{b_j^i\}_{j=1}^k\,$ is independent
we must have $n_1=\cdots =n_k=0$. By definition we then have that $
\{b_j^i\}_{j=1}^\infty $ is an independent subset of $H^i$ and\textbf{\ }$
H^i\cong \bigoplus_{j=1}^\infty \langle b_j^i\rangle _{\mathbb{Z}}\cong
\bigoplus_{j=1}^\infty \mathbb{Z}$. And if we define the $k\times k$ matrix $
M_{i-1}=\left( a_{rs}\right) $ with integer entries by the condition that
for each $r=1,...,k$ we have 
\[
a_{r1}b_1^i+a_{r2}b_2^i+\cdots +a_{rk}b_k^i=b_r^{i-1}\text{,} 
\]
then the entries of $M_{i-1}$ are uniquely determined since $
\{b_j^i\}_{j=1}^\infty $ is independent. The map $\beta _{i-1}^i$ may then
be said to be represented by $\left( M_{i-1}\right) ^T\times id$. That $
M_{i-1}$ is invertible follows from the fact that there are rationals $
q_{rs} $ satisfying 
\[
\left( 
\begin{array}{c}
q_{11}\text{ }\cdots \text{ }q_{1k} \\ 
\vdots \\ 
q_{k1}\text{ }\cdots \text{ }q_{kk}
\end{array}
\right) \left( 
\begin{array}{c}
b_1^{i-1} \\ 
\vdots \\ 
b_k^{i-1}
\end{array}
\right) =\left( 
\begin{array}{c}
b_1^i \\ 
\vdots \\ 
b_k^i
\end{array}
\right) 
\]
by the maximality of $M^{\prime }$ in $A$.\hfill$\square $

Thus, given a $B$-sequence $\{H^i,\beta _i^j\}$ of a countable subgroup of
the reals $H$, where for each $i\in \mathbb{N}$ $H^i\cong
\bigoplus_{j=1}^\kappa \mathbb{Z\,}$ and the bonding map $\beta _i^{i+1}$ is
represented by $\left( M_i\right) ^T$ or $\left( M_i\right) ^T\times id$,
there is the Pontryagin dual inverse sequence $\{G_i,\widehat{\beta _i^j}\}$
, where for each $i\in \mathbb{N}$,  $\;G_i\cong \mathbf{T}^\kappa $ and $
\widehat{\beta _i^{i+1}}$ is represented by $M_i\,$ or $M_i\times id,$ see,
e.g., \cite{K}. Each of the bonding maps $\widehat{\beta _i^{i+1}}$ will be
epimorphic since the corresponding matrix $M_i\,$ is invertible. Hence, the
following is well-defined.

\begin{definition}
\label{B}\emph{We define the }$B$-dual \emph{of the }$B$\emph{-sequence of
the countable subgroup of the reals }$H$\emph{\ }$=\{h_1=0,h_2,...\}$\emph{\
with the bases }$\left\{ b_j^i\right\} $\emph{\ for the }$H^i$\emph{\ of
cardinality }$\kappa $\emph{\ to be the }$\kappa $\emph{-solenoid }$\sum_{
\overline{M}}$\emph{\ which is the inverse limit of the dual inverse
sequence }$\left\{ G_i,\alpha _i^j\right\} $\emph{, where }$G_i=\mathbf{T}
^\kappa $\emph{\ for all }$i\in \mathbb{N}$\emph{\ and the maps }$\alpha
_i^{i+1}$\emph{\ are represented by the transposes of the maps }$\beta
_i^{i+1}$\emph{\ as given with respect to the bases }$\left\{ b_j^i\right\}
_{j=1}^\kappa \,$\emph{\ and }$\left\{ b_j^{i+1}\right\} _{j=1}^\kappa $
\emph{.}
\end{definition}

\begin{theorem}
\label{flowdual}If $f$ is an orbit of an almost periodic flow $\phi $ on $
\left( X,d\right) $ with group of exponents $\mathfrak{M}_f$ $=\{h_1=0,h_2,...\}$
and if $B=\{b_1,b_2,...\}$ is a maximal independent subset of $\mathfrak{M}_f$
and if $\sum_{\overline{M}}$ is the $B$-dual of the $B$-sequence of $\mathfrak{M}
_f$ with the bases $\left\{ b_j^i\right\} $ for the $H^i$ as above, then $
\phi \stackrel{equiv}{\approx }\Phi _{\overline{M}}^{\mathbf{\beta }}$,
where $\mathbf{\beta }=\left( b_1,b_2,...\right) $.
\end{theorem}

\noindent \textbf{Proof: }Since $B$ is independent, $\mathbf{\beta }$ is
irrational and each trajectory of $\Phi _{\overline{M}}^{\mathbf{\beta }}$
is dense in $\sum_{\overline{M}}$. As $\sum_{\overline{M}}$ is compact and
distances are preserved by time $t$ maps of $\Phi _{\overline{M}}^{\mathbf{
\beta }}$, the orbits of $\Phi _{\overline{M}}^{\mathbf{\beta }}$ are almost
periodic, see, e.g., \cite{NS}, V, 8.12, and so we need only show that the group
of exponents $\mathfrak{M}_g$ for the orbit $g$ of $e_{\overline{M}}$ in the
flow $\Phi _{\overline{M}}^{\mathbf{\beta }}$ is the same as $\mathfrak{M}_f$,
see Theorem \ref{class}. In light of Theorems \ref{exp1} and \ref
{exp2}, we need to then show that 
\[
\left[ \{\pi _{\overline{M}}\left( t_i\mathbf{\beta }\right) \}\rightarrow
e_{\overline{M}}\right] \Leftrightarrow \left[ \{\pi \left( \mu t_i\right)
\}\rightarrow \pi \left( 0\right) \text{ in }S^1\text{ for all }\mu \in 
\mathfrak{M}_f\right] \text{.}
\]

Let $\{H^i,\beta _i^j\}$ be the $B$-sequence for $\mathfrak{M}_f$ and choose
bases $\{b_j^i\}_{j=1}^\kappa $ for the groups $H^i$ with $b_j^1\stackrel{\text{def}
}{=}b_j$ for all $j$. This at the same time determines the representation of
the $\beta _i^{i+1}$ by matrices $\left( M_i\right) ^T$. Since the bonding
map $\widehat{\beta _i^{i+1}}$ of $\sum_{\overline{M}}$ is represented by $
M_i$ (or $M_i\times id$\thinspace ), we have: 
\[
\pi _{\overline{M}}\left( t\mathbf{\beta }\right) =\left( \langle \pi \left(
tb_1^1\right) ,\pi \left( tb_2^1\right) ,...\rangle ,\langle \pi \left(
tb_1^2\right) ,\pi \left( tb_2^2\right) ,...\rangle ,...\right) \in
\prod_{i=1}^\infty \mathbf{T}^\kappa \text{.} 
\]
Thus, $\cup _i\{b_j^i\}_{j=1}^\kappa $ generates $\mathfrak{M}_g$ as a subgroup
of $\left( \mathbb{R},+\right) $. And since $\cup _i\{b_j^i\}_{j=1}^\kappa $
generates $\mathfrak{M}_f$ as a subgroup of $\left( \mathbb{R},+\right) $, we have
what we need.\hfill$\square $

Notice that this implies that 
\[
\mathfrak{M}_f=\mathfrak{M}_g\cong \widehat{\sum\nolimits_{\overline{M}}}\cong 
\breve{H}^1\left( \sum\nolimits_{\overline{M}}\right) =\breve{H}^1(\overline{
g\left( \mathbb{R}\right) })\cong \breve{H}^1(\overline{f\left( \mathbb{R}\right) }
). 
\]

	The following corollary of these results gives a specific
form for the inverse limit representation of a compact connected abelian 
topological group. That some inverse limit representation is possible is
guaranteed by a theorem of Pontryagin, \cite{P}, Theorem 68.
\begin{corollary}
\label{group}Any metric compact connected abelian topological group is a $
\kappa $-solenoid for some $\kappa \leq \infty $.
\end{corollary}

\noindent \textbf{Proof}: This follows from the above and the converse in
Nemytskii's Theorem \cite{NS}, V, 8.16. \hfill $\square $

We are now in a position to give our intrinsic classification of minimal
almost periodic flows. The notion of equivalence is generally too strict for
classification, and two flows $\phi $ and $\psi $ on $X$ and $Y$
respectively are said to be \emph{topologically equivalent }(denoted $\phi 
\stackrel{top}{\approx }\psi $) when there is a homeomorphism $
h:X\rightarrow Y\,$ which maps orbits of $\phi $ onto orbits of $\psi $ in
such a way that the orientation of orbits is preserved. This is the notion
of equivalence discussed in \cite{AM}, Chapter 6 for the program of Birkhoff. In
general, equivalence and topological equivalence are different, but in \cite
{LF} we showed that two minimal almost periodic flows are topologically
equivalent if and only if they are equivalent, and any two equivalent flows
have, up to a non-zero multiplicative factor, the same exponent group, see
Corollary \ref{equivexp}. And now that we know that any minimal almost
periodic flow is equivalent with a linear flow on some $\kappa -$solenoid,
we are led to the following conclusion.

\begin{theorem}
\label{classif}If $\phi $ and $\psi $ are minimal almost periodic flows with
exponent groups $\mathfrak{M}$ and $\mathfrak{N}$ respectively, then

\begin{equation*}
\left[ \phi \overset{top}{\approx }\psi \right] \Leftrightarrow \left[ \mathfrak{
M}=a\mathfrak{N}\text{ for some non--zero }a\right] .
\end{equation*}
\end{theorem}

\begin{theorem}
\label{flowdom}If $f=\phi _x$ is the orbit of a flow $\phi :\mathbb{R\times }
X\rightarrow X$ and if $\overline{f\left( \mathbb{R}\right) }$ is compact, then
there is a map $h_f:\overline{f\left( \mathbb{R}\right) }$ $\stackrel{onto}{
\rightarrow }\sum_{\overline{M}}$ with $id\times h_f:\phi \mid _{\mathbb{
R\times }\overline{f\left( \mathbb{R}\right) }}\stackrel{sc}{\succeq }\Phi _{
\overline{M}}^\omega $, where $\sum_{\overline{M}}$ is the $B$--dual of $
\mathcal{E}_f$ and where $B=\left\{ \omega _i\right\} $ is a maximal
independent subset of $\mathcal{E}_f$ and where $\omega =\left( \omega
_1,\omega _2,...\right) $.
\end{theorem}

\noindent \textbf{Proof}: We have the $B$-sequence of $\mathcal{E}_f$\ $
\left\{ H^i,\beta _i^j\right\} $\ with the chosen basis for $H^i$ $\;\left\{
b_j^i\right\} $,$\,$ where $b_j^1=\omega _j$ for all $j$. We then define $
h_f:\overline{f\left( \mathbb{R}\right) }$ $\rightarrow \prod_{\ell =1}^\infty 
\mathbf{T}^\kappa $ as follows, where the maps $f_b:\overline{f\left( \mathbb{R}
\right) }\,\rightarrow S^1$ for $b\in \mathcal{E}_f$ are as in the proof of
Theorem \ref{embedexp}: 
\[
x=\lim \left\{ f\left( t_i\right) \right\} \mapsto \left( \langle f_{\omega
_1}\left( x\right) ,f_{\omega _1}\left( x\right) ,...\rangle ,\langle
f_{b_1^2}\left( x\right) ,f_{b_1^2}\left( x\right) ,...\rangle ,...\right)
\in \prod_{\ell =1}^\infty \mathbf{T}^\kappa , 
\]
where $\kappa $ is the cardinality of $B$. The limits $\lim\limits_i\left\{
\pi \left( b_k^\ell t_i\right) \right\} =f_{b_k^\ell }\left( x\right) $
exist for all $k,\ell $ as shown in Theorem \ref{embedexp}. Also, for each $
i\in \mathbb{N}$ we have that 
\[
\left( \langle \pi \left( \omega _1t_i\right) ,\pi \left( \omega
_2t_i\right) ,...\rangle ,\langle \pi \left( b_1^2t_i\right) ,\pi \left(
b_2^2t_i\right) ,...\rangle ,...\right) \in \sum\nolimits_{\overline{M}} 
\]
by the definition of $\sum_{\overline{M}}$. And so $h_f$ $\left( \overline{
f\left( \mathbb{R}\right) }\right) \subset \sum_{\overline{M}}$ since $\sum_{
\overline{M}}$ is a closed subset of $\prod_{\ell =1}^\infty \mathbf{T}
^\kappa $. Thus, our map $h_f$ is well-defined and it is continuous since
each $f_{b_k^\ell }$ is continuous. We now show that it is surjective. Since 
$\Phi _{\overline{M}}^\omega \,$ is an irrational flow, the orbit of $e_{
\overline{M}}$ is dense, and so given any point $\sigma \in \sum_{\overline{M
}}$ there is a sequence $\{x_i\}\rightarrow \sigma $, where 
\[
x_i=\left( \langle \pi \left( \omega _1t_i\right) ,\pi \left( \omega
_2t_i\right) ,...\rangle ,\langle \pi \left( b_1^2t_i\right) ,\pi \left(
b_2^2t_i\right) ,...\rangle ,...\right) =\Phi _{\overline{M}}^\omega \left(
t_i,e_{\overline{M}}\right) 
\]
and $t_i\in \mathbb{R}$. Since $\overline{f\left( \mathbb{R}\right) }$ is compact,
there is a subsequence $\left\{ t_{i_k}\right\} \,\,$ of $\left\{
t_i\right\} \,$ which is an $f$-sequence. Then $h_f\left( \lim \left\{
f\left( t_{i_k}\right) \right\} \right) =\sigma $ as needed. The
commutativity of the diagram as in Definition \ref{semicon} follows from the
following: for $y=\lim \left\{ f\left( t_i\right) \right\} $, we have 
\begin{eqnarray*}
h_f\circ \phi \left( t,y\right) &=&h_f\left( \lim\limits_i\left\{ f\left(
t+t_i\right) \right\} \right) =\left( \langle \lim_i\left\{ \pi \left(
\omega _1t\right) +\pi \left( \omega _2t_i\right) \right\} ,...\rangle
,...\right) \\
&=&\left( \left\langle \pi \left( \omega _1t\right) ,...\right\rangle
,...\right) +\left( \left\langle f_{\omega _1}\left( y\right)
,...\right\rangle ,...\right) \\
&=&\pi _{\overline{M}}\left( t\omega \right) +h_f\left( y\right) =\Phi _{
\overline{M}}^\omega \left( t,h_f\left( y\right) \right) .
\end{eqnarray*}
\hfill$\square $

We know that if there is an $h$ with $id_{\mathbb{R}}\times h:\overline{f\left( 
\mathbb{R}\right) }\stackrel{sc}{\succeq }\Phi _{\overline{N}}^{\omega ^{\prime
}}$ for some irrational linear flow $\Phi _{\overline{N}}^{\omega ^{\prime }}
$ then $\mathcal{E}_f$ contains the exponent group of $\Phi _{\overline{N}
}^{\omega ^{\prime }}$, see Lemma \ref{dominate}. And so the above
provides a semiconjugacy onto a linear flow which is maximal with respect to
its exponent group.

\section{The Exponent Group Generalizes the Irrational Rotation Number}

\begin{definition}
\label{suspdef}\emph{Let }$h:X\rightarrow X$\emph{\ be a homeomorphism.
Define the equivalence relation }$\approx _h$\emph{\ on }$\mathbb{R\times }X$
\emph{\ by } 
\[
\left[ \left( s,x\right) \approx _h\left( t,y\right) \right] \Leftrightarrow
\left[ \text{\emph{there is an} }n\in \mathbb{Z}\text{\emph{ with } }t=s+n\text{\emph{
and} }y=h^{-n}\left( x\right) \right] 
\]
\emph{and let }$\mathcal{S}_h=\left( \mathbb{R\times }X\right) /\approx _h$
\emph{. The }$\approx _h$\emph{\ class of }$\left( s,x\right) $\emph{\ will
be denoted }$\left[ s,x\right] _h$\emph{. \\ There is then the flow } 
\[
\sigma _h:\mathbb{R\times }\mathcal{S}_h\rightarrow \mathcal{S}_h\text{ ; }
\left( t,\left[ s,x\right] _h\right) \mapsto \left[ t+s,x\right] _h\text{.}
\]
\emph{We refer to both }$\mathcal{S}_h$\emph{ and } $\sigma _h$\emph{ as
the }suspension of $h$.
\end{definition}

The above formulation of suspension may be found, for example, in \cite{A}, but
the notion originated with Smale \cite{Sm}. We now proceed to describe the group
of exponents of the suspension of an orientation preserving homeomorphism of 
$S^1$ with an irrational rotation number, but to do so we will need to relate
suspensions with flows on the torus.

\begin{definition}
For $\theta \in \mathbb{R}$, $R_\theta :S^1\rightarrow S^1$ is the translation given by
\[ 
x\longmapsto x+\pi \left( \theta \right).
\]
\end{definition}

We can then equate the suspension of $R_\theta$ with $\Phi ^{\left( \theta ,1\right) }$
as follows.
\begin{lemma}
\label{mu}
\[
\mu_{R_\theta} : \mathcal{S}_{R_\theta} \longrightarrow \mathbf{T}^2
\]
\[
\left[ s, x \right] _{R_\theta}\longmapsto \left\langle x+\pi \left( s\theta \right),
\pi \left(s \right) \right\rangle 
\]
is a well-defined homeomorphism and
\[
id\times \mu _{R_\theta }:\sigma_{R_\theta }\stackrel{equiv}{\approx }
\Phi ^{\left( \theta ,1\right) } .
\]
\end{lemma}
\noindent See, e.g., \cite{KH}, Proposition 2.2.2 for a proof (ignoring the extraneous
$x_n$).
 
Given a map $f:S^1 \rightarrow S^1$ we may lift the map $f\circ \pi :\mathbb{
R\rightarrow }S^1$ with a map $F:\mathbb{R\rightarrow R}$ satisfying $f\circ
\pi =\pi \circ F,$ which is uniquely determined if we require $F\left(
0\right) \,\,$ to be in $[0,1)$. We assume hereafter that such a choice of $F
$ is made and we refer to $F$ as the lift of $f$. And when $f$ is an
orientation preserving homeomorphism, $F(x+1)=F(x)+1$, see \cite{KH}, Proposition
11.1.1 and so $F(x+k)=F(x)+k$ for $k\in \mathbb{Z}$. And since $F^\ell \,$ is
a lift of $f^\ell $ for $\ell \in \mathbb{Z}$ and $f^\ell $ is an orientation
preserving homeomorphism when $f$ is, we have that for $\left( k,\ell
\right) \in \mathbb{Z}^2$ $F^\ell \left( x+k\right) =F^\ell \left( x\right) +k.$
 Also, $f$ is said to be \emph{monotone} (or \emph{strictly monotone},
etc.) when its lift $F\,$ is monotone (or strictly monotone, etc.) \cite{KH}.
Thus, an orientation preserving homeomorphism of $S^1$ is strictly monotone
increasing.
\newline In \cite{KH} 11.1.2 the rotation number is defined as an 
element of $S^1$. For our purposes it is more convenient to define it 
as a real number. Since we have a unique lift associated with any map, this
rotation number will be well--defined.

\begin{definition}
\emph{Let } $f:S^1\rightarrow S^1$ \emph{ be an orientation-preserving 
homeomorphism with the lift } $F$, \emph{ then the } rotation number \emph{ of } $f$ 
\emph{ is }

\begin{equation*}
\lim_{\left| n\right| \rightarrow \infty }\frac{1}{n}\left( F^{n}\left(
x\right) -x\right) . 
\end{equation*}
\end{definition}
See \cite{KH} 11.1.1.

\begin{definition}
A map $g:N\rightarrow N$ is a \emph{factor} of $
f:M\rightarrow M$ if there exists a surjective map $
h:M\rightarrow N$  such that $h\circ f=g\circ h$. The map $h$
 is then called a \emph{semiconjugacy} \cite {KH}, 2.3.2. If $h$
is a homeomorphism it is said to provide a \emph{conjugacy} and $f$
and $g$ are said to be \emph{conjugate}.
\end{definition}

\begin{theorem}
\label{poincare}\ (Poincar\'{e} Classification Theorem) Let $
f:S^1\rightarrow S^1$ be an orientation-preserving homeomorphism with
irrational rotation number $\theta .$
\end{theorem}

\begin{enumerate}
\item  If $f$ is transitive, then $f$ is conjugate to $
R_\theta .$

\item  If $f$ is not transitive, then the map $R_\theta $ is a factor of $f$
with a non-invertible monotone map $h:S^1\rightarrow S^1$ providing the
semiconjugacy \cite{KH}, 11.2.7.
\end{enumerate}

So in either case there is a surjective monotone map $h$ with $R_\theta
\circ h=h\circ f$.

\begin{lemma}
\label{dom'}If $f:S^1\rightarrow S^1$ is an orientation-preserving
homeomorphism with rotation number $\theta \in \left( \mathbb{R}-\mathbb{Q}\right)
\;$and if $\lambda $ is the $\sigma _f$-orbit of $x\in \mathcal{S}_f$, then $
\mathcal{E}_\lambda \supset \left\langle \theta,1 \right\rangle _{\mathbb{Z}}$.
\end{lemma}

\noindent \textbf{Proof}: By the above theorem we have a monotone map $
h:S^1\rightarrow S^1$ with $R_\theta \circ h=h\circ f$. We then have the map 
\[
\mathcal{S}\left( h\right) :\mathcal{S}_f\rightarrow \mathcal{S}_{R_\theta
};\;\left[ s,x\right] _f\mapsto \left[ s,h\left( x\right) \right] _{R_\theta
}\text{.} 
\]
This map is well-defined since $\left[ s,x\right] _f=\left[ s^{\prime
},x^{\prime }\right] _f$ implies that $s^{\prime }=s+n$ and $x^{\prime
}=f^{-n}\left( x\right) $ for some $n\in \mathbb{Z}$ and: 
\begin{eqnarray*}
\left( \left[ s+n,f^{-n}\left( x\right) \right] _f\right) \stackrel{\mathcal{
S}\left( h\right) }{\mapsto }\left[ s+n,h\left( f^{-n}\left( x\right)
\right) \right] _{R_\theta } &=&\left[ s+n,\left( R_\theta \right)
^{-n}\left( h\left( x\right) \right) \right] _{R_\theta } \\
&=&\left[ s,h\left( x\right) \right] _{R_\theta }=\mathcal{S}\left( h\right)
\left( \left[ s,x\right] _f\right) .
\end{eqnarray*}
We then have the commutative diagram: 
\[
\begin{array}{ccc}
\mathbb{R\times }\mathcal{S}_f & 
\begin{array}{c}
\sigma _f \\ 
\rightarrow
\end{array}
& \mathcal{S}_f \\ 
\begin{array}{cc}
\left( id,\mathcal{S}\left( h \right) \right) & \downarrow
\end{array}
&  & 
\begin{array}{cc}
\downarrow & \mathcal{S}\left( h\right)
\end{array}
\\ 
\mathbb{R\times }\mathcal{S}_{R_\theta } & 
\begin{array}{c}
\sigma _{R_\theta } \\ 
\rightarrow
\end{array}
& \mathcal{S}_{R_\theta }
\end{array}
. 
\]
But we also have the map $\mu _{R_\theta }:\mathcal{S}_{R_\theta
}\rightarrow \mathbf{T}^2$ with $id\times \mu _{R_\theta }:\sigma_{R_\theta }
\stackrel{equiv}{\approx }\Phi ^{\left( \theta ,1\right) }$ (see Lemma \ref{mu}). 
And so 
\[
id\times \left( \mu _{R_\theta }\circ \mathcal{S}\left( h\right) \right)
:\sigma _f\stackrel{sc}{\succeq }\Phi ^{\left( \theta ,1\right) }, 
\]
and by Lemma \ref{dominate} we have $\mathcal{E}_\lambda \supset
\left\langle \theta ,1\right\rangle _{\mathbb{Z}}$ since this is the group of
exponents for each orbit of $\Phi ^{\left( \theta ,1\right) }$.\hfill$\square $

To get results on the other inclusion we first introduce some terminology.

\begin{definition}
\emph{Given a sequence }$\{\xi _i\}$\emph{\ of points in }$S^1$\emph{, we
define }$\{x_i\}\subset \mathbb{R}\,$\emph{\ to be a }lift\emph{\ of }$\{\xi
_i\}$\emph{\ if }$\pi \left( x_i\right) =\xi _i$\emph{\ for all }$i\in \mathbb{N
}$\emph{.}
\end{definition}

\begin{definition}
\emph{A sequence }$\{\xi _i\}$\emph{\ of points in }$S^1$\emph{\ is said to
converge to }$\xi $\emph{\ from below (or above) }$[$\emph{denoted }$
\uparrow $\emph{\ }$\left( \text{\emph{or}}\downarrow \right) ]$\emph{\ if
there is a lift of the sequence }$\{\xi _i\}$\emph{\ to a sequence }$
\{x_i\}\subset \mathbb{R}\,$\emph{\ with }$\{x_i\}\uparrow \left( \text{\emph{or
}}\downarrow \right) $\emph{\ }$x\,$\emph{\ in }$\mathbb{R}$\emph{, where }$\pi
\left( x\right) =\xi $\emph{. We also use }$\uparrow $\emph{\ }$\left(
\downarrow \right) $\emph{\ to indicate that this sort of convergence takes
place.}
\end{definition}

\begin{definition}
\emph{We define a stable orbit }$f=\phi _x:\mathbb{R}\rightarrow X$\emph{\ of a
flow }$\phi $\emph{\ to be }amphiperiodic\emph{\ if the following
implication holds: for some set of generators }$\mathcal{G}$\emph{\ of }$
\mathcal{E}_f$\emph{\ } 
\[
\left[ \{\pi \left( \alpha t_i\right) \}\text{ }\uparrow \text{\emph{ or } }
\downarrow \text{\emph{ for all } }\alpha \in \mathcal{G} \right] \Rightarrow
\left[ \{t_i\}\text{\emph{ is an } }f \text{\emph{--sequence}}\right].
\]
\end{definition}

Notice that requiring $\{\pi \left( \alpha t_i\right) \}$ $\uparrow $ or $
\downarrow $ is stricter than requiring that $\{\pi \left( \alpha t_i\right)
\}$ converges, and so any almost periodic orbit is also amphiperiodic.

\begin{definition}
\emph{If }$f$\emph{\ is an amphiperiodic orbit, we define }$\theta \in 
\mathcal{E}_f$\emph{\ to be }regular\emph{\ if the following implication
holds for some set of generators }$\mathcal{G}$\emph{\ of }$\mathcal{E}_f$
\emph{\ } 
\[
\left( [\{\pi \left( \alpha t_i\right) \}\text{ }\uparrow \text{\emph{or} }
\downarrow \text{\emph{for all} }\alpha \in \mathcal{G}-\left\langle \theta
\right\rangle _{\mathbb{Z}}]\text{ \emph{ and} }[\{\pi \left( \theta t_i\right)
\}\text{ \emph{ converges}}]\right) 
\]
\[
\Rightarrow \left( \{t_i\}\;\text{\emph{is an}\ }f\text{\emph{--sequence}}
\right) ;
\]
\emph{otherwise, }$\theta $\emph{\ is }singular.
\end{definition}

It then follows that all exponents of an almost periodic orbit are regular.

We now quote a standard result we shall need.

\begin{theorem}
Let $f:S^1\rightarrow S^1$ be an orientation-preserving
homeomorphism with rotation number $\theta \in \left( \mathbb{R}-\mathbb{Q}
\right) $. Then there is a unique minimal set $M_f$ for the
dynamical system 
\[
\mathbb{Z}\times S^1\rightarrow S^1;\;\left( n,x\right) \mapsto f^n\left(
x\right) ,
\]
which is either $S^1$ or a perfect and nowhere dense subset
(i.e., a Cantor set), see, e.g., \cite{KH}, 11.2.5.
\end{theorem}

By the Poincar\'{e} classification, $M_f$ is $S^1$ when $f\,\,$ is transitive
and a Cantor set otherwise. This Cantor set has as its complement a
countable union of pairwise disjoint open intervals, each of which gets
mapped by $h$ to a single point. The endpoints of these open intervals get
mapped to the same point as the open interval by $h$, while the remainder of 
$M_f$ is mapped one-to-one by $h$ onto its image, see, e.g., \cite{KH}, p. 398.
Then the set $\mathcal{M}_f=\left\{ \sigma _f\left( t,\left[ 0,m\right]
_f\right) \mid t\in \mathbb{R}\text{ and }m\in M_f\right\} $ is a minimal set
of the flow $\sigma _f$, see, e.g., \cite{Schw}. Any aperiodic $C^1$ flow on $
\mathbf{T}^2$ is topologically equivalent to the suspension of an
orientation-preserving circle diffeomorphism with irrational rotation
number, see, e.g., \cite{KH} 14.2.3, 0.3 and 11.1.4, and so any minimal set
occurring in an aperiodic $C^1$ flow on $\mathbf{T}^2$ that is a proper
subset of $\mathbf{T}^2$ is homeomorphic with some such $\mathcal{M}_f$.

\begin{theorem}
\label{amphi}If $f:S^1\rightarrow S^1$ is a non-transitive,
orientation-preserving homeomorphism with rotation number $\theta \in \left( 
\mathbb{R}-\mathbb{Q}\right) $ and if $\lambda $ is the orbit of a point $x=\left[
s,m\right] _f\in \mathcal{M}_f$, then $\mathcal{E}_\lambda =\left\langle
\theta ,1\right\rangle _{\mathbb{Z}}$ and $\lambda \,$ is amphiperiodic and $1$
is a regular exponent while $\theta $ is a singular exponent.
\end{theorem}

\noindent \textbf{Proof}: Suppose then that for a given sequence $
\{t_i\}\subset \mathbb{R}$ we have that $\{\pi \left( \theta t_i\right) \}$ $
\uparrow $ or $\downarrow $ and that $\{\pi \left( t_i\right) \}$ converges.
We shall show that $\{t_i\}$ is a $\lambda $-sequence. Since we know that $
\left\{ \theta ,1\right\} \subset \mathcal{E}_\lambda $, this will
demonstrate that $\lambda $ is amphiperiodic. As before, we have the
monotone map $h:S^1\rightarrow S^1$ with $h\circ f=R_\theta \circ h$ and $
id\times \mu _{R_\theta }:\sigma _{R_\theta }\stackrel{equiv}{\approx }\Phi
^{\left( \theta ,1\right) }$. Then with $\mu _{R_\theta }\circ \mathcal{S}
\left( h\right) \left( x\right) =\left\langle \xi ,\pi \left( u\right)
\right\rangle \in \mathbf{T}^2$ and $\mathfrak{t}$ translation by $\left\langle
-\xi ,-\pi \left( u\right) \right\rangle $ and $g$ defined to be the map $
\mathfrak{t}\circ \mu _{R_\theta }\circ \mathcal{S}\left( h\right) $, we have
that $id\times g:\sigma _f\stackrel{sc}{\succeq }\Phi ^{\left( \theta
,1\right) }$ since $\Phi ^{\left( \theta ,1\right) }$ is translation
invariant and $g\circ \lambda \left( t\right) =\left\langle \pi \left(
\theta t\right) ,\pi \left( t\right) \right\rangle \in \mathbf{T}^2$. Thus,
our hypothesis on $\{t_i\}$ guarantees that $\lim\limits_ig\circ \lambda
\left( t_i\right) =\lim\limits_i\left\langle \pi \left( \theta t_i\right)
,\pi \left( t_i\right) \right\rangle $ exists and is equal to say $
\left\langle \chi ,\pi \left( \tau \right) \right\rangle $. Since $h$ is an
onto monotone function, the set $h^{-1}\left( \xi +\chi \right) $ is either 
(\textbf{A1}) a point or\textbf{\ }(\textbf{A2) }a closed interval.

\begin{description}
\item[(A1)]  $h^{-1}\left( \xi +\chi \right) =\pi \left( r\right) $.
\end{description}

In this case we claim that $\left\{ \lambda \left( t_i\right) \right\}
\rightarrow \left[ \tau +u,\pi \left( r\right) \right] _f$. Let 
\[
U_\varepsilon \stackrel{\text{def}}{=}\left\{ \left[ \ell +\tau +u,\pi \left(
\left( r-\varepsilon ,r+\varepsilon \right) \right) \right] _f\in \mathcal{S}
_f\mid \ell \in \left( -\varepsilon ,\varepsilon \right) \right\} \text{.} 
\]
Then $\left\{ U_\varepsilon \mid \varepsilon >0\right\} $ forms a local base
at $\left[ \tau +u,\pi \left( r\right) \right] _f$. And for any given $
\varepsilon >0$ the intervals $\pi \left( (r-\varepsilon ,r)\right) $ and $
\pi \left( \left( r,r+\varepsilon \right) \right) $ are not mapped by $h$ to
single points, and so their image under the monotone map $h$ contains an open
interval. Thus, the image of $U_\varepsilon $ under $g$ contains a
neighborhood of $g\left( \left[ \tau +u,\pi \left( r\right) \right]
_f\right) =\left\langle \chi ,\pi \left( \tau \right) \right\rangle $ and so
contains the tail of the sequence $\left\{ \left\langle \pi \left( \theta
t_i\right) ,\pi \left( t_i\right) \right\rangle \right\} $. Thus, the tail
of the sequence $\left\{ \lambda \left( t_i\right) \right\} $ is contained
in  $\widetilde{U_\varepsilon }\stackrel{\text{def}}{=}\left( g\right)
^{-1}\left( g\left( U_\varepsilon \right) \right) $. While $\widetilde{U_\varepsilon }$ 
properly contains 
$U_\varepsilon $ if one of the endpoints $\pi \left( r-\varepsilon \right) $
or $\pi \left( r+\varepsilon \right) $ is in the closure of a complementary
open interval, in this case $\widetilde{U_\varepsilon }-U_\varepsilon $ does
not contain points of arbitrarily high index $i$ from the set $\left\{
\lambda \left( t_i\right) \right\} $ since $\widetilde{U_\varepsilon }
-U_\varepsilon $ gets mapped to a line segment (or two) in $\mathbf{T}^2$
which is (are) a positive distance from $\left\langle \chi ,\pi \left( \tau
\right) \right\rangle $, and so all but finitely many terms of the sequence $
\left\{ \lambda \left( t_i\right) \right\} $ are contained in $U_\varepsilon 
$. Notice that in this case we did not use the full strength of the
hypothesis that $\{\pi \left( \theta t_i\right) \}$ $\uparrow $ or $
\downarrow $ , we only needed that $\{\pi \left( \theta t_i\right) \}$
converges -- a fact we use later.

\begin{description}
\item[(A2)]  $h^{-1}\left( \xi +\chi \right) =\pi \left( \left[ \mu ,\nu
\right] \right) $.
\end{description}

Suppose without loss of generality that $\{\pi \left( \theta t_i\right)
\}\uparrow \chi $ and that that $h$ is monotone increasing, 
then we claim that $\left\{ \lambda \left( t_i\right)
\right\} \rightarrow \left[ \tau +u,\pi \left( \mu \right) \right] _f$.
Since, for $\varepsilon >0$, $h$ does not collapse the interval $\pi \left(
(\mu -\varepsilon ,\mu ]\right) $ to a point, $h$ maps each such interval
onto a set containing an open interval. For each $\varepsilon >0$ define

\[
U_\varepsilon =\left\{ \left[ \ell +\tau +u,\pi \left( (\mu -\varepsilon
,\mu ]\right) \right] _f\in \mathcal{S}_f\mid \ell \in \left( -\varepsilon
,+\varepsilon \right) \right\} . 
\]
Then $g\left(U_\varepsilon \right) \,$ contains an open subset of the torus and $
\left\langle \chi ,\pi \left( \tau \right) \right\rangle $, and, since $
\{\pi \left( \theta t_i\right) \}$ $\uparrow \chi $, $g\left( U_\varepsilon
\right) \,$ contains the tail of the sequence $\left\{ \left\langle \pi
\left( \theta t_i\right) ,\pi \left( t_i\right) \right\rangle \right\} $.
Then $U_\varepsilon $ contains a tail of the sequence since $g^{-1}\left(
g\left( U_\varepsilon \right) \right) -U_\varepsilon $ contains only
finitely members of the sequence $\left\{ \lambda \left( t_i\right) \right\} 
$ as above. Notice that if the tail of the sequence $\left\{ \left\langle
\pi \left( \theta t_i\right) ,\pi \left( t_i\right) \right\rangle \right\} $
also contained points with $\chi +\pi \left( \theta t_i\right) ``>"\chi $,
then the inverse image of such points would not be in $U_\varepsilon $ since
these inverse images would be ``on the other side'' of the complementary
disk. And since each neighborhood of $\left[ \tau +u,\pi \left( \mu \right)
\right] _f$ contains some $U_\varepsilon $, we have that $\left\{ \lambda
\left( t_i\right) \right\} \rightarrow \left[ \tau +u,\pi \left( \mu \right)
\right] _f$. In this case we made full use of our hypotheses.

This shows then that $\lambda $ is amphiperiodic. It only remains to show $
\mathcal{E}_\lambda \subseteq \left\langle \theta ,1\right\rangle _{\mathbb{Z}} ,$
 and the proof of this bears some resemblance to the proof of a related fact for almost
periodic orbits (see, e.g., \cite{LZ}, Chapter 2.2.3, p. 42-3). Let $\gamma \notin
\left\langle \theta ,1\right\rangle _{\mathbb{Z}}$. We have two cases: (\textbf{B1})
$\{\gamma ,1,\theta \}$ is
rationally independent and (\textbf{B2}) $\{\gamma ,1,\theta \}$ is
rationally dependent.

\begin{description}
\item[(B1)]  $\{\gamma ,1,\theta \}$ is rationally independent.
\end{description}

Let $\left\langle \zeta ,\pi \left( 0\right) \right\rangle \,\in \mathbf{T}
^2\,$ be a point for which $\left( g\right) ^{-1}\left( \left\langle \zeta
,\pi \left( 0\right) \right\rangle \right) $ is a single point, which will
be the case whenever $\zeta $ is the image under $h\left( x\right) -\xi $ of
a point $x$ of the Cantor set $M_f$ which in not an endpoint. By Kronecker's
Theorem we can construct a sequence $\{t_i\}\,$ of real numbers with $
\left\{ \pi \left( t_i\right) \right\} \rightarrow \pi \left( 0\right) $ and 
$\left\{ \pi \left( \theta t_i\right) \right\} \rightarrow \zeta $ and such
that $d_1\left( \pi \left( t_i\gamma \right) ,\pi \left( 0\right) \right)
<\frac 1i$ for odd $i$ and $d_1\left( \pi \left( t_i\gamma \right) ,\pi
\left( \frac 12\right) \right) <\frac 1i$ for even $i$. Then we have that: 
\[
\left\{ \left\langle \pi \left( \theta t_i\right) ,\pi \left( t_i\right)
\right\rangle \right\} _i\rightarrow \left\langle \zeta ,\pi \left( 0\right)
\right\rangle . 
\]
By our choice of $\left\langle \zeta ,\pi \left( 0\right) \right\rangle $
and Case (\textbf{A1}), we have that $\{t_i\}$ is a $\lambda $-sequence (we
have that both $\left\{ \pi \left( \theta t_i\right) \right\} $ and $\left\{
\pi \left( t_i\right) \right\} $ converge and this is all we need in this
case). Therefore, $\gamma \notin \mathcal{E}_\lambda $ since $\left\{ \pi
\left( t_i\gamma \right) \right\} $ does not converge.

\begin{description}
\item[(B2)]  $\{\gamma ,1,\theta \}$ is rationally dependent. 
\end{description}

There are thus integers
$\ell _1,\ell _2$ and $\ell _3$ not all 0 such that $
\ell _1\cdot 1+\ell _2\theta +\ell _3\gamma =0$. And since 
$\{1,\theta \}$ is rationally independent, $\ell _3$ cannot be 0.
Then we have a representation with the least positive $q$: 
\[
\gamma =\frac pq+\frac rq\theta \text{, for integers }p,q\text{ and }r.
\]
Let $\alpha _1\in \pi ^{-1}\left( 0\right) $ and $\alpha _2\in \pi
^{-1}\left( \zeta \right) $ and let 

\begin{equation*}
\alpha _{3}=\frac{p}{q}\alpha _{1}+\frac{r}{q}\alpha _{2} ,
\end{equation*}
where $\left\langle \zeta ,\pi \left(
0\right) \right\rangle $ is as in \textbf{(B1)}. Then if for integers $\ell
_1,\ell _2$ and $\ell _3$ we have that $\ell _1\cdot 1+\ell _2\theta +\ell
_3\gamma =0$, we also have 
\[
\left( \frac{\ell _1}{\ell _3}+\frac pq\right) +\left( \frac{\ell _2}{\ell _3
}+\frac rq\right) \theta =0\text{,} 
\]
and so the rational independence of $\{1,\theta \}$ yields that 
\begin{equation*}
\frac{p}{q}=-\frac{\ell _{1}}{\ell _{3}}\text{ and }\frac{r}{q}=-\frac{\ell
_{2}}{\ell _{3}} .
\end{equation*}

And so 

\begin{equation*}
\pi \left( \ell _{1}\alpha _{1}+\ell _{2}\alpha _{2}+\ell _{3}\alpha
_{3}\right) =\pi \left( \ell _{1}\alpha _{1}+\ell _{2}\alpha _{2}+\ell
_{3}\left( -\frac{\ell _{1}}{\ell _{3}}\alpha _{1}+-\frac{\ell _{2}}{\ell
_{3}}\alpha _{2}\right) \right) =\pi \left( 0\right) .
\end{equation*}

Thus, the
conditions of the strong version of Kronecker's Theorem (see, e.g., \cite{LZ},
Chapter 3.1, p. 37) are met and so we may find for any $\delta >0$ a real number 
$t$ satisfying: 
\[
\;d_1\left( \pi \left( t\right) ,\pi \left( \alpha _1\right) \right) <\delta
,\text{ }d_1\left( \pi \left( \theta t\right) ,\pi \left( \alpha _2\right)
\right) <\delta \;\text{and }d_1\left( \pi \left( \gamma t\right) ,\pi
\left( \alpha _3\right) \right) <\delta \text{.} 
\]
Now let 
\begin{equation*}
\alpha _{3}^{\prime }=\frac{1}{q}+\frac{p}{q}\alpha _{1}+\frac{r}{q}\alpha_{2} .
\end{equation*}

Then we have that $\pi \left( \alpha
_3^{\prime }\right) \neq \pi \left( \alpha _3\right) $ since $q$ cannot be $
1 $; otherwise, $\gamma $ would be in $\left\langle \theta ,1\right\rangle _{
\mathbb{Z}}$. And if for integers $\ell _1,\ell _2$ and $\ell _3$ we have that $
\ell _1\cdot 1+\ell _2\theta +\ell _3\gamma =0$, as before we have 

\begin{equation*}
\frac{p}{q}=-\frac{\ell _{1}}{\ell _{3}}\text{ and }\frac{r}{q}=-\frac{\ell
_{2}}{\ell _{3}}
\end{equation*}

and by our choice of $q$ we must
have that $q$ divides $\ell _3$. And so:

\begin{eqnarray*}
\pi \left( \ell _{1}\alpha _{1}+\ell _{2}\alpha _{2}+\ell _{3}\alpha
_{3}^{\prime }\right)  &=&\pi \left( \ell _{1}\alpha _{1}+\ell _{2}\alpha
_{2}+\ell _{3}\left( \frac{1}{q}-\frac{\ell _{1}}{\ell _{3}}\alpha _{1}+-%
\frac{\ell _{2}}{\ell _{3}}\alpha _{2}\right) \right)  \\
&=&\pi \left( \frac{\ell _{3}}{q}\right) =\pi \left( 0\right) .
\end{eqnarray*}

Again we have the conditions for the strong version of Kronecker's Theorem
and so we may find for any $\delta >0$ a real number $t$ satisfying: 
\[
\;d_1\left( \pi \left( t\right) ,\pi \left( \alpha _1\right) \right) <\delta
,\text{ }d_1\left( \pi \left( \theta t\right) ,\pi \left( \alpha _2\right)
\right) <\delta \;\text{and }d_1\left( \pi \left( \gamma t\right) ,\pi
\left( \alpha _3^{\prime }\right) \right) <\delta \text{.} 
\]
Thus, we may construct a sequence $\{t_i\}$ with $\left\{ \pi \left(
t_i\right) \right\} \rightarrow \pi \left( \alpha _1\right) =\pi \left(
0\right) $ and $\left\{ \pi \left( \theta t_i\right) \right\} \rightarrow
\pi \left( \alpha _2\right) =\zeta $ and such that 
\[
d_1\left( \pi \left( t_i\gamma \right) ,\pi \left( \alpha _3\right) \right)
<\frac 1i\text{ for odd }i 
\]
and 
\[
d_1\left( \pi \left( t_i\gamma \right) ,\pi \left( \alpha _3^{\prime
}\right) \right) <\frac 1i\text{ for even }i. 
\]
Then 
\[
\left\{ \left\langle \pi \left( \theta t_i\right) ,\pi \left( t_i\right)
\right\rangle \right\} _i\rightarrow \left\langle \zeta ,\pi \left( 0\right)
\right\rangle 
\]
and by our choice of $\left\langle \zeta ,\pi \left( 0\right) \right\rangle $
and Case (\textbf{A1}), we have that $\{t_i\}$ is a $\lambda $-sequence.
Therefore, $\gamma \notin \mathcal{E}_\lambda $ since $\left\{ \pi \left(
t_i\gamma \right) \right\} $ does not converge.\hfill$\square $

Closer examination of the map $g$ reveals that it is nothing other than $
h_\lambda \,$, where $h_\lambda \,$is as in Theorem \ref{flowdom}. The
complement of $\mathcal{M}_f$ in $\mathcal{S}_f$ consists of ``blown up''
trajectories of $\Phi ^{\left( 1,\theta \right) }$, see \cite{KH}, p. 398, each
component of which is homeomorphic to an open disk. To see that the
component of $\mathcal{S}_f-\mathcal{M}_f$ containing $\left\{ \left[
0,x\right] _f\in \mathcal{S}_f\mid x\in J\right\} $, where $J$ is one of the
countable pairwise disjoint intervals forming a component of $S^1-M_f$, is
homeomorphic with $D=\mathbb{R}\times J$, we construct the one-to-one map $
\delta :D\rightarrow \left( \mathcal{S}_f-\mathcal{M}_f\right) ;$
\[
\left( t,x\right) \mapsto \left[ t,x\right] _f
\]
(that $\delta $ is one-to-one follows from the fact that $\sigma _f$ has no
periodic orbits). Brouwer's Theorem on the Invariance of Domain yields that $
\delta $ is an open map and thus a homeomorphism onto its image, see also
\cite{Fok}, p. 47. To investigate $\check{H}^1\left( \mathcal{M}_f\right) $ we
examine a portion of the long exact reduced singular homology sequence for
the pair $\left( \mathcal{S}_f,\mathcal{S}_f-\mathcal{M}_f\right) $, see,
e.g., \cite{Br}, p.184-5: 
\[
H_1\left( \mathcal{S}_f-\mathcal{M}_f\right) \rightarrow H_1\left( \mathcal{S
}_f\right) \rightarrow H_1\left( \mathcal{S}_f,\mathcal{S}_f-\mathcal{M}
_f\right) \rightarrow \tilde{H}_0\left( \mathcal{S}_f-\mathcal{M}_f\right)
\rightarrow \tilde{H}_0\left( \mathcal{S}_f\right) \text{.}
\]
Since each component of $\mathcal{S}_f-\mathcal{M}_f$ is homeomorphic to an
open disk, $H_1\left( \mathcal{S}_f-\mathcal{M}_f\right) =0$. And so $\tilde{
H}_0\left( \mathcal{S}_f-\mathcal{M}_f\right) \cong \oplus _{i=1}^{\kappa -1}
\mathbb{Z}$, where as above $\kappa $ represents the number of orbits blown up
to form $f$, and so there is an exact sequence of groups 
\[
0\rightarrow \mathbb{Z}^2\rightarrow H_1\left( \mathcal{S}_f,\mathcal{S}_f-
\mathcal{M}_f\right) \rightarrow \oplus _{i=1}^{\kappa -1}\mathbb{Z\rightarrow }
0 .
\]
By Poincar\'{e}-Lefschetz Duality (see, e.g., \cite{Br}, VI, 8.4), we have an
isomorphism $\check{H}^1\left( \mathcal{M}_f\right) \cong H_1\left( \mathcal{
S}_f,\mathcal{S}_f-\mathcal{M}_f\right) $. Thus, for sufficiently large
$\kappa ,$ $\left[ \mathcal{M}_f;S^1\right] \cong \check{H}^1\left( \mathcal{M}_f\right) $ 
 has torsion-free rank greater than two and so must properly
contain the image of the exponent group $\iota \left( \mathcal{E}_\lambda
\right) ,$ which has  torsion-free rank two.

Now we proceed to generalize exponent groups to self--homeomorphisms.

\begin{definition}
\emph{Let }$h:X\rightarrow X$\emph{\ be a homeomorphism. Then for }$x\in X$
\emph{\ we define the group }$\mathcal{E}_{\left( x,h\right) }$\emph{\ to be
the group }$\mathcal{E}_f$\emph{\ for the orbit }$f$\emph{\ of the point }$
\left[ 0,x\right] _h$\emph{\ in the suspension }$\mathcal{\sigma }_h$ .
\end{definition}

Notice that $\mathcal{E}_{\left( y,h\right) }$ is the same for all $y$ in
the $h$-orbit of $x$. And if $M$ is a minimal set of $h$,
all points of $M$ have the same exponent group. By our calculation in
Theorem \ref{amphi}, this exponent group generalizes the rotation number in
some sense. 

\begin{example}
Let $h: U \rightarrow \mathbb{R}^2 $ be an elementary 
twist mapping of a domain of $ U \subset \mathbb{R}^2 $ as described in 
\cite{AM} 8.3.5. 
\end{example}
Such a map occurs naturally as the return map to a cross section of a Hamiltonian flow 
restricted to a 3--dimensional energy surface. While $h$ itself may not be a 
self--homeomorphism of $U$, we may restrict $h$ to the interior $V$ of an invariant circle
so that the restricted map $h'$ will be a self--homeomorphism of $V$, see \cite{AM} 8.3.6.
For any invariant circle in $V$ as described in part (i) of Moser's Theorem \cite{AM} 8.3.6,
we can calculate the rotation number or the exponent group.
For some such $h'$ there will also be minimal Cantor sets on which $\sigma_{h'}$  
is equivalent to a linear flow $\Phi _{\overline{N}
}^{1}$ on a 1--sonlenoid $\sum_{\overline{N}}$, see \cite{AM} pp. 583-586. And so if, for 
example, $\sum_{\overline{N}}$ is the dyadic solenoid, the exponent group for all points
in the corresponding Cantor set of $V$ will be the dyadic rationals. Of course, if $h'$ actually
is the return map of a flow $\phi ,$ this does not imply that the exponent group of the 
$\phi$--orbits of points
in the homeomorphic copy of $\sum_{\overline{N}}$ will be the corresponding subgroup of the 
rationals: $\phi$ would be topologically equivalent to $\sigma_{h'}$  
but not necessarily equivalent.

We end with a question.
 
\begin{description}
\item[Question] Is there a non-singular compact minimal set $M$ for a
flow whose exponent group is $\left\{ 0\right\} $?
\end{description}

\end{document}